\documentclass[12pt]{amsart}
\usepackage{amsmath,amscd,amssymb,amsfonts}
\newcommand{\h}{\hbox}
\newcommand{\q}{\quad}
\newcommand{\nin}{\par\noindent}
\newcommand{\bs}{\par\bigskip}
\newcommand{\ms}{\par\medskip}
\newcommand{\sk}{\par\smallskip}
\newcommand{\bsn}{\par\bigskip\noindent}
\newcommand{\msn}{\par\medskip\noindent}
\newcommand{\skn}{\par\smallskip\noindent}

\newcommand{\mopl}{\h{$\bigoplus$}}
\newcommand{\mcap}{\h{$\bigcap$}}

\newcommand{\ssb}{\raise.15ex\h{${\scriptscriptstyle\bullet}$}}
\newcommand{\ssc}{\,\raise.15ex\h{${\scriptstyle\circ}$}\,}
\newcommand{\CC}{{\mathcal C}}
\newcommand{\D}{{\mathcal D}}
\newcommand{\DD}{{\mathbb D}}
\newcommand{\HH}{{\mathcal H}}
\newcommand{\M}{{\mathcal M}}

\newcommand{\R}{{\mathcal R}}
\newcommand{\X}{{\mathcal X}}
\newcommand{\OO}{{\mathcal O}}
\newcommand{\C}{{\mathbb C}}
\newcommand{\DDD}{{\mathbb D}}
\newcommand{\LL}{{}\,\overline{\!L}{}}

\newcommand{\PP}{{\mathbb P}}
\newcommand{\Q}{{\mathbb Q}}
\newcommand{\RRR}{{\mathbb R}}

\newcommand{\U}{{\mathcal U}}
\newcommand{\Z}{{\mathbb Z}}
\newcommand{\RR}{\widetilde{\mathcal R}}
\newcommand{\FFc}{{}\,\overline{\!F^c}{}}
\newcommand{\MM}{{}\,\overline{\!M}{}}
\newcommand{\MMc}{{}\,\overline{\!M^c}{}}
\newcommand{\Db}{{\mathfrak D}{\mathfrak b}}
\newcommand{\XX}{{}\,\overline{\!X}{}}
\newcommand{\dd}{\partial}
\newcommand{\ddd}{{\rm d}}
\newcommand{\al}{\alpha}
\newcommand{\la}{\lambda}

\newcommand{\Om}{\Omega}

\newcommand{\Gr}{{\rm Gr}}
\newcommand{\bl}{\bigl}
\newcommand{\br}{\bigr}
\newcommand{\into}{\hookrightarrow}

\newcommand{\simto}{\buildrel\sim\over\too}
\newcommand{\too}{\longrightarrow}
\newcommand{\ges}{\geqslant}

\begin{document}
\title{Mixed Hodge modules and mixed twistor modules}
\author[M. Saito]{Morihiko Saito}
\begin{abstract} We explain some fundamental differences between the theories of mixed Hodge modules and mixed twistor modules (including the difference in weight system on the nearby cycle functor) which do not seem to be clarified explicitly in the literature.
\end{abstract}
\maketitle
\centerline{\bf Introduction}
\bsn
In the introduction of the first version of \cite{Sab4}, it is stated that there is a {\it fully faithful\,} functor from the category of mixed Hodge modules ${\rm MHM}(X)$ (see \cite{mhm}) to the category of mixed twistor modules ${\rm MTM}(X)$ (see \cite{Mo2}).
If $X$ is a point, then this implies that there is a fully faithful functor from the category of mixed (complex) Hodge structures to the category of mixed twistor structures in the sense of Simpson \cite{Si} (see also \cite[2.1]{Sab2}).
The ``full faithfulness" implies that the {\it Hodge numbers} can be recovered from the image of a mixed Hodge structure under this functor.
However, it is shown by Simpson \cite{Si} that this can be done only by using a $\C^*$-action.
This ``full faithfulness" does not seem to be attained even by replacing ${\rm MTM}(X)$ with the category of {\it integrable} mixed twistor modules ${\rm MTM}^{\rm int}(X)$ as in the second version of \cite{Sab4}, if the integrability condition is taken in a weak sense as in \cite[7.1]{Sab2}, \cite[2.8--9]{Sab4}. Indeed, it does not seem enough to assume that a mixed twistor module {\it admits} an action of $z^2\dd_z$ (or $\lambda^2\dd_{\lambda}$ in \cite{Mo2}) satisfying certain properties as in \cite{Sab2}, \cite{Sab4}, but we would have to {\it fix} such an action in order to capture the Hodge filtration $F$ as in \cite{Si}. In fact, without it, ${\rm MTM}^{\rm int}(X)$ is still a {\it full subcategory} of ${\rm MTM}(X)$, and ${\rm MTM}^{\rm int}(X)$ would not be called naturally a {\it subcategory} of ${\rm MTM}(X)$ as in the introduction of \cite{Sab4} unless one can {\it choose} naturally a good action of $z^2\dd_z$ for each object in this subcategory, although this does not seem quite easy (since the situation is entirely different from the ``integrability" of a connection), see Remark~(1.6) below for more detailed explanations. Anyway it seems rather difficult to say that the natural functor from ${\rm MTM}^{\rm int}(X)$ to ${\rm MTM}(X)$ is a {\it forgetful functor} forgetting the action of $z^2\dd_z$.
\sk
As an example, consider pure Hodge structures of rank $2$ and type
$$\{(p,-p),(-p,p)\}\q\h{for}\,\,\,\,p\in\Z_{>0},
\leqno(1)$$
(see \cite{th2}), which are denoted by $H_{(p)}$. It is quite unclear how one can find a difference between their images for $p\in\Z_{>0}$ in the category of twistor structures of weight 0 without using a $\C^*$-action (which is essentially equivalent to a $z\dd_z$-action) as in \cite{Si}.
Recall that the category of twistor structures of weight $k$ consists of vector bundles $E$ on $\PP^1$ which are (non-canonically) isomorphic to direct sums of copies of $\OO_{\PP^1}(k)$, see \cite{Si}, \cite[2.1]{Sab2}.
\sk
Here it is unclear whether a polarization of a pure Hodge structure using the Weil operator is {\it correctly} compared in \cite[Lemma~3.46]{Mo1} with a polarization of the corresponding pure twistor structure endowed with a $\C^*$-action or a $z\dd_z$-action.
Notice that there are two ways of {\it sign convention} for polarizations of Hodge structures used by Deligne and Griffiths, where the difference comes from the place of the Weil operator $C$, and produces a {\it certain difference} in the sign of polarizations of Hodge structures depending on the {\it dimensions of strict supports} (more precisely, see \cite[1.4.7]{ypg}).
Indeed, Deligne's convention is used in mixed Hodge modules (see for instance \cite[Remark before Theorem~3.20]{mhm}), although Griffiths' one is mainly used in Hodge theory recently (including papers of Kashiwara and Kawai), see for instance \cite[Section 5.5]{St}.
It does not seem very clear which convention is used in \cite{Mo1}, \cite{Mo2}, \cite{Sab2}, since the notation is rather complicated.
\sk
It is also usually unnoticed that there is a fundamental difference in the {\it weight system} on the nearby cycle functor for Hodge modules and twistor modules, and the weight of a pure twistor module does not not change under non-characteristic restriction, see Section~2 below.
\sk
I would like to thank T.~Mochizuki for very important information \cite{Mo3}, \cite{Mo4} giving some explanations about the correspondence between Hodge and twistor modules in \cite[13.5]{Mo2} and also about \cite[Lemma 3.7.9]{Sab2}.
\bs\bs
\centerline{\bf 1. ``Integrable $\R$-triples" or ``$\RR$-triples"}
\bsn
{\bf 1.1.~$\RR$-modules.} Let $X$ be a complex manifold. Set $\X=X\times\C$. By definition $\R_{\X}$ is the subalgebra of $\D_{\X}$ generated by $z\,p_1^*\xi$ for $\xi\in\Theta_X$ over $\OO_{\X}$, where $z$ is the coordinate of $\C$, $p_1:\X\to X$ is the first projection, and $\Theta_X$ is the sheaf of vector fields. This $\R_{\X}$ is contained in the subalgebra $\RR_{\X}$ of $\D_{\X}$ generated by $z^2\dd_z$ over $\R_{\X}$.
\sk
An $\R$-triple consists of $(\M',\M'',C)$ where the $\M',\M''$ are $\R_{\X}$-modules and $C$ is a certain pairing between $\M'|_{X\times S}$ and $\sigma^*\M''|_{X\times S}$ where $S:=\{z\in\C\mid|z|=1\}$ and $\sigma$ is induced by the involution of $\C^*$ defined by $z\mapsto-1/\overline{z}$.
It is called {\it integrable} if the $\R_{\X}$-module structure of $\M',\M''$ is {\it liftable} to an $\RR_{\X}$-module structure so that the pairing $C$ satisfies a certain compatibility condition for the action of $z^2\dd_z$ (see \cite[7.1]{Sab2}).
\sk
It seems much better to use the terminology ``$\RR$-triple", rather than ``integrable $\R$-triple" here, since the $\RR$-module structure is not uniquely determined by the underlying $\R$-module structure as is shown in (1.2) below. This seems rather misleading for non-specialists.
\msn
{\bf 1.2.~Case $X=pt$.} The above problem becomes clearer in the case $X=pt$. This case was studied by C.~Simpson \cite{Si}, who showed that $\C^*$-action is needed to capture the Hodge filtration. Here we have
$$\R_{\X}=\OO_{\C},\q\RR_{\X}=\OO_{\C}\langle z^2\dd_z\rangle,$$
with $z$ the coordinate of $\C$.
For simplicity, assume furthermore that $\M',\M''$ are finite free $\OO_{\C}$-modules and the action of $z^2\dd_z$ comes from that of $z\dd_z$, that is, the corresponding connection has a logarithmic pole.
(Note that a $\C^*$-action is essentially equivalent to a $z\dd_z$-action.)
The eigenvalues of the {\it residue} of the logarithmic connection should be closely related to the complex numbers $p$, $q$ with
$$\Gr_F^pi_1^*\M'\ne 0,\,\,\,\Gr_F^qi_1^*\M''\ne 0\q\h{(counted with multiplicities)},$$
where $i_1:\{1\}\into\C$ is the inclusion, see \cite{Sab2}.
(Note that these are closely related to the ambiguities of the integrable structure in \cite[Remark~2.9]{Sab4}.)
\sk
Indeed, the twistor structure associated with a complex Hodge structure $(H;F',F'')$ can be defined by using
$$\mopl_{p\in\Z}\,\,\overline{\!F'^pH}\otimes z^{-p},\q\mopl_{q\in\Z}\,F''^qH\otimes z^{-q}.
\leqno(1.2.1)$$
More precisely, we have to use the {\it dual} of the first term, since the {\it first} term is {\it contravariant} in the case of twistor modules, and a {\it polarization} of complex Hodge structure is {\it not} used here as is noted in an earlier version of this note, see \cite[2.1]{Sab2}, and also (A.2--4) below. Indeed, if we do not use the dual, then we would get a twistor structure of weight 0, since the dual Hodge stricture $H^{\vee}$ is isomorphic to the conjugate Hodge structure $\,\overline{\!H}$ up to the Tate twist $(w)$ with $w$ the weight of the Hodge structure.
\sk
The above argument shows that there is a {\it canonical representative} of $\R$-modules (up to a non-canonical isomorphism) using the action of $z^2\dd_z$ in the case of twistor modules coming from variation of complex Hodge structures.
\msn
{\bf Remark~1.3.} The graded $\C[z]$ modules in (1.2.1) are quite similar to Brieskorn lattices of Gauss-Manin systems associated with {\it weighted homogeneous} polynomials having isolated singularities, see \cite{Sab3}, \cite{bl}, \cite{ScSt}, etc.) In the general ``integrable" case where $z^2\dd_z$-actions are simply assumed, one may have a situation similar to the Brieskorn lattices of certain {\it non-weighted-homogeneous} polynomials.
Here one has to use the graded-quotients of the $V$-filtration along $z=0$. This may be related to the rescaling limit argument in \cite{Sab4} in the case $X=pt$, since the latter seems to be related to the theory of asymptotic Hodge structures by Varchenko \cite{Va} where one takes the {\it leading terms} of asymptotic integrals, see also \cite{ScSt}.
(Note that the irregular case with respect to $\dd_z$ is also possible if we consider $z^2\dd_z$-actions.)
\msn
{\bf 1.4.~Description of the Hodge filtration $F$.} The argument in (1.2) implies that the Hodge filtration $F$ {\it cannot be captured} unless one {\it fixes} an $\RR$-module structure for $\M',\M''$. Indeed, the eigenvalues of the residues depend heavily on the lifting as an $\RR$-module, see also the example in the introduction which gives {\it various liftings as $\RR$-module structures} for various $p\in\Z_{>0}$. (This may be related with \cite[2.8--9]{Sab4}.) So it would be better to use ``$\RR$-triple" rather than ``integrable $\R$-triple" in order to avoid any possible confusions.
(Indeed, ``integrable $\R$-triple" may strongly suggest that an $\RR$-module structure is not fixed although it admits such a structure, see also Remark~(1.6) below.)
\sk
Related to the above problem, it does not seem very clear how to interpret, for instance, an ``integrable morphism of the underlying filtered $\R$-triples" in \cite[p.~188, $\ell$.~5]{Mo2}. It is noted in \cite[7.1.c]{Sab2} that a morphism of integrable triples is called ``integrable" if it commutes with ``some representatives of the $\bar{\dd}_z$-actions" (with $\bar{\dd}_z:=z^2\dd_z$).
However, it does not seem quite clear whether some ``equivalence class" is fixed in \cite[p.~188, $\ell$.~5]{Mo2} or not.
(Two actions of $\bar{\dd}_z$ on $(\M',\M'')$ are called ``equivalent" in \cite[7.1.c]{Sab2}, if their difference is given by $\lambda z$ on $\M'$ and $\bar{\lambda}z$ on $\M''$ for some $\lambda\in\C$.)
This is closely related to \cite[2.8-9]{Sab4} in the {\it irreducible} (or {\it simple}) case.
Here it does not seem very clear whether {\it direct sums} of integrable triples can be well-defined under the above definition of integrable morphisms, since the ambiguity $\lambda$ might depend on direct factors.
This problem seems to be closely related to the example in the introduction, since the Hodge structures there are direct sums of two {\it complex Hodge structures}, see \cite{Si}, \cite[2.1.d]{Sab2}.
\sk
As a conclusion, it does not seem quite easy to show that the integrability condition in the weak sense explained above together with the real structure is sufficient to get a fully faithful functor from the category of mixed Hodge modules.
Note that the full faithfulness implies that, for an object in the essential image, the corresponding mixed Hodge module could be determined up to a canonical isomorphism; in particular, the Hodge numbers could be determined uniquely in the case $X=pt$.
\msn
{\bf 1.5.~Relation with Kashiwara's conjecture.}
It seems {\it highly desirable} to construct a certain category of twistor modules on any variety in such a way that this category contains the category of Hodge modules {\it as a full subcategory} (or at least as a {\it subcategory}) and moreover the underlying $\D$-modules of its objets contain {\it any} irreducible holonomic $\D$-modules. However, one cannot cover all the irreducible local systems on smooth complex varieties if he imposes the above integrability condition (in the weak sense as explained above) on twistor deformations.
In the case of rank 1 local systems on smooth projective varieties, for instance, the integrability condition implies that the twistor deformation (which is slightly different from Simpson's $\C^*$-action) is constant so that the Higgs field vanishes, and hence the local system is {\it unitary}, see \cite{Sab2} and also \cite{tdef} (here it is not very clear whether one can get a correct definition of twistor deformations of local systems in the non-compact case without assuming the extendability to a compactification, see \cite[Remark~2.5(iii)]{tdef}).
So it would be better to avoid the same notation MTM for the one with integrability condition imposed (as in someone's talk in Kyoto).
Perhaps ``(irregular) Hodge-twistor modules" might be more suitable for these objects, since they seem to be quite close to ``irregular Hodge modules", see \cite{Sab4}.
\msn
{\bf Remark~1.6.} In the introduction of (the first version of) \cite{Sab4}, it is stated as follows: {\it ``The subcategory ${\rm MTM}^{\rm int}(X)$ of integrable objects and morphisms plays an important intermediate role in what follows. The $\R_{\X}$-modules underlying the objects in this subcategory are equipped with a compatible action of $z^2\dd_z$ and the pairing is supposed to be compatible with it. The morphisms in this subcategory are those morphisms in ${\rm MTM}(X)$ which are compatible with the $z^2\dd_z$."}
Here the most nontrivial point is the following question: Of which category is ${\rm MTM}^{\rm int}(X)$ considered as a ``{\it subcategory}"? In view of the above expression, this seems to be regarded as a ``subcategory" of ${\rm MTM}(X)$, and there does not seem to be no other choices. Indeed, one does not seem to be talking about the category whose objects are objects of ${\rm MTM}(X)$ endowed with a good action of $z^2\dd_z$, that is, the category of $W$-filtered $\RR_{\X}$-triples satisfying certain good conditions (the latter category will be denoted by ${\rm MTM}^{\sim}(X)$ in these notes). In fact, ``this subcategory" would be replaced by ``this category" in such a case. Moreover the last sentence of the quoted phrases seems to describe the way in which the author {\it shrinks} the {\it groups of morphisms} of ``this subcategory".
However, there seems to be {\it some non-trivial difference} between the category ${\rm MTM}^{\sim}(X)$ and the situation mentioned in the above quoted sentences. It does not seem quite easy to {\it realize} the above situation without solving some {\it set-theoretic problem} as is explained below.
\sk
Indeed, in order to realize the above sentence, one would have to ``give", or rather ``choose", an action of $z^2\dd_z$ (among all the possible choices) for each object in this ``subcategory" of ${\rm MTM}(X)$.
This seems to be equivalent to choosing an {\it objectwise section} (forgetting about morphisms) on the image of the forgetful functor ${\rm MTM}^{\sim}(X)\to{\rm MTM}(X)$. The image of the latter coincides with the {\it objects} of ${\rm MTM}^{\rm int}(X)$. However, this forgetful functor is never injective except for certain special cases, since there is an {\it ambiguity} of choice of the action of $z^2\dd_z$. This means the the same object of ${\rm MTM}(X)$ may underlie many different objects of ${\rm MTM}^{\sim}(X)$. It implies that the natural functor from ${\rm MTM}^{\sim}(X)$ to ${\rm MTM}(X)$ {\it never factors through} ${\rm MTM}^{\rm int}(X)$, and the {\it groups of morphisms of $\,{\rm MTM}^{\rm int}(X)$ do depend} on the choice of the objectwise section, that is, on the choice of the action of $z^2\dd_z$. In particular, ${\rm MTM}^{\rm int}(X)$ would be called rather ``{\it a} subcategory" instead of ``{\it the} subcategory" unless the choice of the $z^2\dd_z$-action is explicitly given.
(In order to understand the situation, the reader may restrict to the case $X=pt$ and consider only pure objects of weight 0 as in the example in the introduction.
It may be also helpful to compare the situation to the case where one identifies a {\it certain full subcategory} of the category ${\bf M}(k[x])$ of $k[x]$-modules with a {\it ``subcategory"} of the category ${\bf M}(k)$ of $k$-vector spaces by {\it choosing} an element of ${\rm End}_k(V)$ as an action of $x$ for each $V\in{\bf M}(k)$, where $k$ is a field and $k[x]$ is the polynomial ring in one variable $x$. Here one gets only a {\it full subcategory} of ${\bf M}(k[x])$ since for each $k$-vector space $V$, there is only {\it one object} in this full subcategory such that its underlying $k$-vector space is $V$. Note that this full subcategory of ${\bf M}(k[x])$ is not necessarily equivalent to the whole category ${\bf M}(k[x])$ if the above ``subcategory" of ${\bf M}(k)$ is badly given.)
\sk
It does not seem very clear, however, how the above choice is made in \cite{Sab4}.
This seems to be rather nontrivial even in the case of mixed twistor structures where $X=pt$.
There does not seem to be a {\it canonical} way to do it in general.
It does not seem even clear whether there is a {\it good} way to do it, although it seems relatively easy to choose a {\it very bad} action of $z^2\dd_z$ in such a way that for any two objects in this subcategory which are isomorphic to each other in ${\rm MTM}(X)$, the actions of $z^2\dd_z$ on these two objects are compatible with {\it some} isomorphism between these in ${\rm MTM}(X)$. (Indeed, choose an object in each isomorphism class of objects of ${\rm MTM}(X)$, choose one isomorphism between the chosen object and any other object in the isomorphism class, and then choose an action of $z^2\dd_z$ for the chosen object in the isomorphism class.
In the case of twistor structures of weight $0$, this would correspond to assigning, for instance, a trivial Hodge structure of type $(0,0)$ to {\it any} twistor structure of weight $0$.)
It does not seem very clear in general how these bad choices are excluded in the above quoted situation, for instance, if we choose the action of $z^2\dd_z$ {\it arbitrarily} among possible choices for each object of ``a subcategory" ${\rm MTM}^{\rm int}(X)$.
\sk
Actually it may be possible to give a {\it very artificial choice} of the action of $z^2\dd_z$ by dividing each isomorphism class of objects of ${\rm MTM}(X)$ into a disjoint union indexed by all the possible actions of $z^2\dd_z$ on the chosen object in this isomorphism class (provided that the associated set-theoretical problem can be resolved).
In this case, however, it seems more appropriate to say that ``a subcategory" ${\rm MTM}^{\rm int}(X)$ of ${\rm MTM}(X)$ is identified with a certain ``full subcategory" of ${\rm MTM}^{\sim}(X)$ which is equivalent to ${\rm MTM}^{\sim}(X)$, by using this ``very artificial choice".
Here the condition: ``which is equivalent to ${\rm MTM}^{\sim}(X)$" disappears if we take a bad choice of the $z^2\dd_z$-action as above.
Anyway, how to choose an action of $z^2\dd_z$ on each object of ``a subcategory" ${\rm MTM}^{\rm int}(X)$ does not seem to be a trivial matter which can be left without giving any details.
\msn
{\bf Remark~1.7.} In the third version of \cite{Sab4}, it is stated as follows: {\it ``The category ${\rm MTM}^{\rm int}(X)$ of integrable objects and morphisms plays an important intermediate role in what follows. The $\R_{\X}$-modules underlying the objects in this category are equipped with a compatible action of $z^2\dd_z$ and the pairing is supposed to be compatible with it. The morphisms in this category are defined as are the morphisms in ${\rm MTM}(X)$, with the supplementary condition that they are compatible with the $z^2\dd_z$-action."} Here it does not seem very clear whether the problem explained in Remark~1.6 above is completely solved by this change. The main problem is that the {\it same object} of ${\rm MTM}(X)$ {\it can} underlie {\it many different objects} of ${\rm MTM}^{\rm int}(X)$, or ${\rm MTM}^{\sim}(X)$ following the notation in Remark~(1.6) above, as is repeatedly explained there. In view of this phenomenon, the expression: ``The $\R_{\X}$-modules underlying the objects in this category are equipped with a compatible action of $z^2\dd_z$" could be rather confusing, since it is not quite clear whether the $z^2\dd_z$-action on the underlying $\R_{\X}$-modules {\it can really depend on each object of ${\rm MTM}^{\rm int}(X)$} (or ${\rm MTM}^{\sim}(X)$ in the notation of Remark~(1.6)). Indeed, it seems quite possible for the reader to understand that the author is simply considering ``the $\R_{\X}$-modules underlying the objects in this category" without caring so much about the object of this category that each $\R_{\X}$-module underlies. It may be better to note, for instance, as follows: ``The $\R_{\X}$-modules underlying each object in this category are equipped with a compatible action of $z^2\dd_z$ depending on each object, and not only on the underlying $\R_{\X}$-modules". (Here note that a pair of $\R_{\X}$-modules is needed for each object of ${\rm MTM}(X)$.)
\sk
The above problem may have some relation to the difference between an ``inverse functor" and a ``quasi-inverse functor" for a functor of categories $F:\CC\to \CC'$. Indeed, the former is a functor $G:\CC'\to\CC$ such that $G\ssc F$ and $F\ssc G$ are the {\it identity} functors on $\CC$ and $\CC'$. This is, however, rather difficult to construct in practice. What we can usually construct is the latter which satisfies the following conditions: There are functorial isomorphisms $\psi_A:G(F(A))\cong A$ and $\phi_B:F(G(B))\cong B$ for $A\in\CC$, $B\in\CC'$.
\bs\bs
\centerline{\bf 2. Weight system of mixed twistor modules}
\bsn
{\bf 2.1.~Weights on nearby cycles.} Let $X$ be a smooth variety.
The weight of the structure sheaf $\OO_X$ is always 0 {\it independently} of the dimension of $X$, according to the authors of \cite{Sab2} and \cite{Mo2}. Here we consider the twistor module corresponding to the {\it constant} twistor deformation of $\OO_X$ over $\PP^1$ or a constant variation of twistor structure of rank 1 and weight 0, see \cite{Sab2}, \cite{Si}, and also Remark~(2.8) below.
The above assertion follows from the definition of twistor modules using the nearby cycle functors along holomorphic functions (see \cite[4.1]{Sab2} and also \cite[p.~9, $\ell$.~10]{Mo2}), which we apply to local coordinates of $X$ inductively. Indeed, according to these, the weight filtration $W$ on the nearby cycle functor $\psi_f$ of a pure twistor module $\M$ of weight $w$ is given by the monodromy filtration shifted by $w$, instead of $w-1$ as in the Hodge module case. Note that the latter shift of weight by $-1$ in the mixed Hodge module case implies that the weight of a pure Hodge module with a strict support is the sum of the dimension of the support and the pointwise weight at a generic point of the support as in the $\ell$-adic case (see also Kashiwara's remark explained below).
\msn
{\bf 2.2.~Kashiwara's remark.} If we define the weight filtration $W$ on the vanishing cycle functor $\varphi_{f,1}$ in the same way as $\psi_f$ (that is, the monodromy filtration $W$ is shifted by the weight $w$), then we get the {\it vanishing} of
$${\rm can}:\Gr^W_k\psi_{f,1}\M\to\Gr^W_k\varphi_{f,1}\M,\q{\rm var}:\Gr^W_k\varphi_{f,1}\M\to(\Gr^W_k\psi_{f,1}\M)(-1),
\leqno(2.2.1)$$
see Kashiwara's remark before \cite[Theorem 3.21 on p.~303]{mhm}. This shift for $\varphi_{f,1}$ is quite reasonable if one considers the case where the support of $\M$ is contained in $f^{-1}(0)$.
\msn
{\bf 2.3.~Weights on direct images.} It turns out that the weights of the direct images of twistor modules under projective morphisms would be given (as is noted in \cite[Theorem 6.1.1]{Sab2}) as follows:
\msn
(A)\q If $f:X\to Y$ is a projective morphism and $\M$ is pure of weight $w$, then $\HH^jf_+\M$ is pure of weight $w+j$.
\ms
This seems to hold, for instance, in the case where
$$\M=\OO_X\,\,\,\h{with}\,\,\,w=0,\,\,\,Y=pt,\,\,\,\h{and}\,\,\,\HH^0f_+\OO_X=\HH^0(a_X)_+\OO_X\,\,\,\h{has weight}\,\,\,0.
\leqno(2.3.1)$$
Here $a_X:X\to pt$ denotes the structure morphism, and we consider the constant twistor deformation of $\OO_X$ over $\PP^1$ as in (2.1). Recall that the direct image functor of twistor modules for the projection $a_X:X\to pt$ is defined by using the relative de Rham complex which is locally the Koszul complex of the $z\dd_{x_i}$ using the bases $z^{-p}\,\ddd x_{i_1}\wedge\cdots\wedge\ddd x_{i_p}$, where the division by $z^p$ corresponds to the shift of the Hodge filtration $F$. We have, for instance, the induced pairing between $H^0(X,\Om_X^{d_X})\otimes z^{-d_X}$ and its complex conjugate, where $d_X:=\dim X$. This seems to define a twistor structure of weight $0$, and not $d_X$, as is explained a remark after (1.2.1). Indeed, there is a canonical pairing between $H^0(X,\Om_X^{d_X})\otimes z^0$ and its complex conjugate, which is {\it constant} on $S^1:=\{z\in\C\mid|z|=1\}$. This would imply the assertion for $H^0(X,\Om_X^{d_X})\otimes z^{-d_X}$ and its complex conjugate, since $z\bar z=1$ on $S^1$.
\sk
More generally, $(p,q)$-forms have a perfect pairing with $(d_X-p,d_X-q)$ forms, and the conjugates of the former are $(q,p)$ forms. This would imply that $H^j(a_X)_+\OO_X$ has weight $j$. Indeed, $H^j(a_X)_*\OO_X$ should be defined by
$$\bl(H^{d_X-j}(X,\C),H^{d_X+j}(X,\C),C\br),$$
(see \cite[1.6.13]{Sab2}), and we have
$$q-(d_X-p)=j\q\h{if}\q p+q=d_X+j.$$
(This fundamental example does not seem to be explained in the literature.)
\sk
In the ``integrable" case (more precisely, in the ``$\RR$-triple case"), the above argument should give the ``Hodge numbers" of the direct image.
\sk
The assertion (A) implies that the weights of twistor modules {\it do not} change under the direct images by closed embeddings. (This is {\it different} from the earlier version of this note.)
\msn
{\bf Remark~2.4.} There is a Tate twist $(-1/2)$ (in the notation of \cite{Sab2}) in the isomorphism between the two nearby cycle functors in \cite[Proposition 4.3.1]{Mo2}, where the usual one is as defined in \cite{Sab2}, and is denoted simply by $\psi$ in this note, and the other one, which is denoted by $\psi^{(1)}$, is defined in \cite{Mo2} by using Beilinson's construction. It may be more natural to use the latter for the inductive definition of twistor modules from the beginning in \cite{Sab2}.
\msn
{\bf Remark~2.5.} In the calculation of the nearby cycles for a special case which was used in the proof of the decomposition theorem (see \cite[Lemma 3.7.9]{Sab2}), we would have
$$\Gr_k^W\psi_f\M=\begin{cases}(i_Z)_+\M|_Z(1/2)&\h{if}\,\,\,k=w-1,\\(i_{Y_1})_+\M|_{Y_1}\oplus(i_{Y_2})_+\M|_{Y_2}&\h{if}\,\,\,k=w,\\(i_Z)_+\M|_Z(-1/2)&\h{if}\,\,\,k=w+1,\end{cases}
\leqno(2.5.1)$$
Here $w$ is the weight of a pure twistor module $\M$, $f^{-1}(0)$ is a union of normally crossing two smooth hypersurfaces $Y_1,Y_2$ in $X$, and $Y_1,Y_2,Z:=Y_1\cap Y_2$ are transversal to a Whitney stratification of $\M$.
\sk
The above Tate twists are not stated in \cite{Sab2} (by the reason that the weights of twistor modules can be changed arbitrarily since there are Tate twists with half-integer values, that is, with any integer weights, {\it according to the author of} \cite{Sab2}).
It turns out that {\it the above Tate twists do not appear at all on the level of $\R$-modules, and they are reflected only on the pairings $C$} (in the manner noted at the end of \cite[Lemma 3.7.9]{Sab2}) according to \cite{Mo4}.
\sk
In the ``integrable" case (or more precisely, in the $\RR$-triple case), however, the above assertion may sound rather strange, since there would be a {\it canonical} representative of $\RR$-modules using the action of $z^2\dd_z$ at least in the variation of complex Hodge structure case as is explained in (1.2). In this case, some more arguments concerning the action of $z^2\dd_z$ may be needed for the proof of \cite[Lemma 3.7.9]{Sab2}. (For instance, the Tate twists in (2.5.1) must be {\it more precise} in the ``integrable" case, see (2.9) below.)
\sk
Even in the non-integrable case, the argument in \cite{Sab2} (and \cite{Mo4}) may be a little bit confusing because of the definition of $P\psi$ in \cite[3.6.14]{Sab2} since the relation with $\M|_Z$, the usual restriction of $\M$ to $Z$, is not stated there (indeed, the assertion for the underlying $\R$-modules is proved for the first time in the first part of \cite[Lemma 3.7.9]{Sab2}). So the last part of \cite[Lemma 3.7.9]{Sab2} may be misstated. Actually $\Gr_1^W\psi_fC$ seems to be compared with the pairing of $\M|_Z$ by using {\it two-variable} Mellin transform in the {\it proof} of \cite[Lemma 3.7.9]{Sab2}. However, the argument is rather difficult to follow for non-experts, since the assertion of \cite[Lemma 3.7.8]{Sab2} states only the information about the orders of poles and some assertion from the proof of \cite[Proposition 3.8.1]{Sab2} (that is, (3.8.2) which is shown only in the constant variation case) is quoted.
\msn
{\bf Remark~2.6.} Let $f:X\to S$ be a projective morphism of smooth varieties with $\dim S=1$. Assume $D:=f^{-1}(0)$ is a reduced divisor with simple normal crossings for $0\in S$. Let $D_i$ be the irreducible components of $D$. Set
$$D_I:=\mcap_{i\in I}\,D_i\q\h{with $\,\,i_I:D_I\into X\,\,$ the natural inclusion.}$$ Let $t$ be a local coordinate of $(S,0)$. Then the {\it co-primitive} part of the graded quotient of weight $-k$ of the weight filtration on $\psi_f\OO_X:=\psi_t(i_f)_+\OO_X$ is described as
$$\bl(\Gr^W_{-k}\psi_f\OO_X\br){}_{\rm copr}=\mopl_{|I|=k+1}\,(i_I)_+\OO_{D_I}(k/2)\q(k\ges 0),
\leqno(2.6.1)$$
and the {\it primitive} part of the graded quotient of weight $k$ as
$$\bl(\Gr^W_{k}\psi_f\OO_X\br){}_{\rm prim}=\mopl_{|I|=k+1}\,(i_I)_+\OO_{D_I}(-k/2)\q(k\ges 0),
\leqno(2.6.2)$$
where the latter is related to the former by $N^k$, see \cite[Proposition 3.8.1]{Sab2}. (Note, however, that $i_J^*$ used there is the usual restriction, and is quite different from the pull-back functor in (2.6.4) below, see also \cite[14.3.3]{Mo2}.)
\sk
Note that the co-primitive part of $\Gr^W_{\ssb}\psi_f\OO_X$ should be related closely to the $\D$-module corresponding to $\C_Y$ by the {\it local invariant cycle theorem} asserting that
$$\aligned\Q_Y[\dim Y]=\,&{\rm Ker}\bl(N:\psi_{f,1}\Q_X\to\psi_{f,1}\Q_X(-1)\br)\\\bl(=\,&{\rm Ker}({\rm can}:\psi_{f,1}\Q_X\to\varphi_{f,1}\Q_X)\br).\endaligned
\leqno(2.6.3)$$
These suggest that the native restriction morphism does not work in twistor theory. More precisely, for a closed embedding $i:X\into Y$ of smooth complex manifolds of codimension $r$, and for a twistor module $\M$ on $Y$ which is {\it non-characteristic} to $X$ (for instance, $X$ is transversal to a Whitney stratification of $\M$), the pull-back functor $i^*$ in the derived category (see \cite[14.3.3]{Mo2}) should satisfy
$$i^*\M\cong\M|_X(r/2)[r],
\leqno(2.6.4)$$
where $M|_X$ denotes the usual restriction as in \cite{Sab2}.
(This may be explained somewhere in the literature, see \cite[Lemma 14.3.26]{Mo2} for the Tate twisted constant variation case.) It is closely related to Remark~(2.7) below.
\sk
Note that (2.6.4) suggests that, for a smooth morphism $f:X\to Y$ of relative dimension $r$, we would have
$$f^*\M\cong f^{\ssb}\M(-r/2)[-r],
\leqno(2.6.5)$$
where $f^{\ssb}\M$ denotes the usual smooth pull-back as in \cite{Sab2}. In the case $Y=pt$, the Tate twist in (2.6.5) seems to be compatible with \cite[13.5]{Mo2}, see (2.8.1) below.
\msn
{\bf Remark~2.7.} Let $X$ be a smooth projective curve, and $P$ be a point of $X$ with $i_P:\{P\}\into X$ the inclusion. In this note, we denote $a_X^*\Q\in D^b{\rm MHM}(X)$ by $\Q_{h,X}$ with $a_X:X\to pt$ the structure morphism (and similarly for $\Q_{h,P}$). There is a canonical morphism (coming from adjunction morphisms)
$$(\Q_{h,X}[1])\to(\Q_{h,P})[1],
\leqno(2.7.1)$$
in the derived category of mixed Hodge modules on $X$, where $(M)$ indicates that $M$ is a Hodge module for $M=\Q_{h,X}[1]$ or $\Q_{h,P}$ (and $(i_P)_*$ before $\Q_{h,P}$ is omitted to simplify the notation). More precisely, the distinguished triangle associated to this morphism is expressed by the following short exact sequence in the abelian category of mixed Hodge modules ${\rm MHM}(X)$:
$$0\to\Q_{h,P}\to(j_U)_!\Q_{h,U}[1]\to\Q_{h,X}[1]\to 0,
\leqno(2.7.2)$$
where $U:=X\setminus\{P\}$ with inclusion $j_U:U\into X$.
\sk
Taking the direct image of (2.7.1) by $a_X:X\to pt$, we then get the canonical isomorphism
$$H^{-1}(a_X)_*(\Q_{h,X}[1])=H^0(X,\Q)\simto H^{-1}(a_X)_*(\Q_{h,P}[1])=\Q.
\leqno(2.7.3)$$
\skn
Here one problem is as follows:
\msn
\rlap{(Q)}\q\q Is there a morphism corresponding to (2.7.1) in the bounded derived category\par\nin\q\q
of mixed twistor modules?
\ms
It does not seem that there is a corresponding morphism in twistor theory if one considers $\OO_X$ (with constant twistor deformation over $\PP^1$) and its usual restriction to $P\in X$. Here one would have to use Beilinson's nearby cycle functor $\psi^{(1)}$ to define the restriction to $P$. This functor is different from the usual one in \cite{Sab2} (that is, $\psi$ in this note) by the Tate twist $(1/2)$ as is explained in Remark~(2.4). These are closely related to (2.6.4).
\sk
Note also that the twistor module $\OO_X(-d_X/2)$ in the notation of \cite{Sab2} is used in \cite[13.5]{Mo2} as the twistor module corresponding to the Hodge module $\Q_{h,X}[d_X]$. Here the notation is rather complicated since the dual filtered $\D$-modules are used, see also Remark~(2.8) below. (These follow by interpreting Mochizuki's comments \cite{Mo3} given {\it very recently}.)
If one uses $\widetilde{\OO}_X:=\OO_X(-d_X/2)$ instead of $\OO_X$ in Remark~(2.6), then one would have to use Beilinson's nearby cycle functor instead of the usual one to get a good formula.
\msn
{\bf Remark~2.8.} If $(\OO_{\X},\OO_{\X},C)$ denotes the twistor module of weight 0 corresponding to the constant twistor deformation of $\OO_X$ over $\PP^1$ as in (2.1) (where $\X:=X\times\C$), then $\OO_X(-d_X/2)$ in Remark~(2.7) above is isomorphic to
$$(z^{d_X}\OO_{\X},\OO_{\X},C)\cong(\OO_{\X},z^{-d_X}\OO_{\X},C)\br)\,\bl(\cong(\OO_{\X},\OO_{\X},C)(-d_X/2)\br)
\leqno(2.8.1)$$
where $z^{d_X}$ in the first term corresponds to the shift of the Hodge filtration $F$ (that is, the Tate twist). Note that the first isomorphism in (2.8.1) would not hold unconditionally in the ``integrable" case (that is, in the $\RR$-triple case), since this would change the Hodge numbers if the representative is taken canonically by using the action of $z^2\dd_z$ as in the last remark in (1.2).
\sk
The above shift of the Hodge filtration comes from the self-duality isomorphism
$$\DD(\OO_X,F)=(\OO_X,F[d_X]).
\leqno(2.8.2)$$
A similar isomorphism holds for any pure Hodge modules with $d_X$ replaced by the weight $w$ in general, and this is used in  \cite[13.5]{Mo2} in an essential way. This may work at least if one uses {\it right} $\D$-modules where the filtration $F$ on $\Om_X^{d_X}$ is shifted by $-d_X$. This shift may imply some shift of the filtration $F$ in \cite[13.5]{Mo2}.
If we use {\it left} $\D$-modules instead, then the filtration $F$ on the dualizing sheaf is shifted by $2d_X$, and this may induce a shift of the filtration $F$ in the duality isomorphism for the direct images by projective morphisms, which is used in an essential way in \cite[13.5]{Mo2}. (Note also that there is a shift of filtration $F$ in the transformation between left and right $\D$-modules, and this shift depends on the ambient dimension. This will also induces a shift of filtration $F$ in the duality isomorphism for the direct images by proper morphisms.)
It is rather complicated for the author to follow the arguments in \cite[13.1.1--2]{Mo2} (see also \cite{Mo4}).
\msn
{\bf 2.9.~Tate twists.} In the ``integrable" case (that is, in the $\RR$-triple case), {\it Tate twists} should be indexed by {\it two} integers $p,q$. For instance, the ``integrable" twistor structure associated to a complex Hodge structure of rank 1 as in (1.2) should have some type $(p,q)$ coming from the Hodge-type of the Hodge structure, and the Tate twists should contain the information about the change of types. These {\it refined Tate twists} may be denoted, for instance, by $(\!(a,b)\!)$, where $-a,\,-b$ give the {\it shift} of type $(p,q)$ (that is, $H(\!(a,b)\!)$ has type $(p-a,q-b)$ if $H$ is an ``integrable" twistor structure of type $(p,q)$). For $\M=(\M',\M'',C)$, these can be defined by
$$\M(\!(a,b)\!):=(z^{-a}\M',z^b\M'',C).
\leqno(2.9.1)$$
They coincide with the Tate twists $\U(-a,b)$ in \cite[2.1.8.1]{Mo2}, and give $((a+b)/2)$ in \cite{Sab2} {\it forgetting} the action of $z^2\dd_z$. In the ``integrable" case, it may be possible to {\it define} that the Tate twist $(m)$ means $(\!(m,m)\!)$ for $m\in\Z$, but this is {\it unclear} for $m\in\tfrac{1}{2}\,\Z\setminus\Z$. (Note that the Tate twist $(m)$ changes the weights by $-2m$.)
\sk
Note that, by considering the argument in (2.3), it may be possible to use also
$$\M(\!(a,b)\!)^t:=(z^{-b}\M',z^a\M'',C),
\leqno(2.9.2)$$
where $\M(\!(a,b)\!)^t=\M(\!(b,a)\!)$. (In this case, (1.2.1) might be modified appropriately.)
\sk
Anyway the Tate twists in (2.5.1), (2.6.1--2), (2.6.4--5) should be replaced by the above {\it refined} Tate twists in the ``integrable" case.
\msn
{\bf Remark~2.10.} As for the definition of twistor modules, it seems necessary to {\it assume} that the pairing $C$ of twistor modules depends {\it holomorphically} on $z\in S^1:=\{z\in\C^*\mid\,|z|=1\}$ in some sense; more precisely, it should depend {\it real analytically} by using an automorphism of $\PP^1$ moving $S^1$ to the real axis $\RRR\cup\{\infty\}\subset\PP^1$. Indeed, we would have to glue {\it holomorphic vector bundles}, and it does not seem easy to replace these by $C^{\infty}$ bundles. As for \cite[Lemma 1.5.3]{Sab2}, it does not seem very clear how ``Grothendieck's Dolbeault lemma" is used for its proof if the definition of ${\mathcal C}_{{\mathcal X}|S}^{\infty,{\rm an}}$ is really as in \cite[Example 0.5.2]{Sab2}; for instance, in the case $X=pt$.
\sk
As to the $E_2$-degeneration argument of the weight spectral sequence, the same remark as in (A.6) below may apply.
\bs\bs
\centerline{\bf Appendix. Some remarks on complex Hodge modules}
\bsn
Recently the theory of {\it complex} Hodge modules is studied by some people. These modules are between real Hodge modules and twistor modules, and the theory seems to be based on some ideas of Kashiwara (see \cite{Ka}, \cite{Sab1}, \cite[(5.1)]{ScVi}). We note here some remarks related to them.
\msn
{\bf A.1.~Complex conjugation.}
Let $X$ be a complex manifold, and $\OO_X$ be the sheaf of holomorphic functions. We denote respectively by $\OO_{\XX}$, $\OO_{X_{\RRR}}$ the sheaf of anti-holomorphic functions and complex-valued real analytic functions (or $C^{\infty}$ functions if one prefers) on the underlying topological space of $X$.
\sk
Let $L$ be a $\C$-local system on $X$. Its complex conjugate is denoted by $\LL$. This is defined by the sheaf $L$ with action of $\C$ given via the {\it complex conjugation}
$$c\mapsto\overline{c}\q(c\in\C).
\leqno({\rm A}.1.1)$$
If $L$ has rank 1 and a local monodromy of $L$ is given by the multiplication by $\la\in\C^*$, then the corresponding local monodromy of $\LL$ is the multiplication by $\overline{\la}$.
\sk
Set
$$M:=\OO_X\otimes_{\C}L.$$
This is a locally free $\OO_X$-module with holomorphic connection $\nabla$. Set
$$\MM:=\OO_{\XX}\otimes_{\C}\LL,$$
which is the complex conjugate of $M$. This is a locally free $\OO_{\XX}$-module with anti-holomorphic connection $\overline{\nabla}$, and is defined by the sheaf $M$ with action of $\OO_{\XX}$ given via the {\it complex conjugation}
$$g\mapsto\overline{g}\q(g\in\OO_{\XX}).
\leqno({\rm A}.1.2)$$
\msn
{\bf A.2.~Variations of complex Hodge structure.}
With the above notation, assume the local system $L$ underlies a polarizable variation of {\it complex} Hodge structure. Here the Hodge filtration $F$ is defined on $M$. We have the {\it opposite Hodge filtration} $F^c$ which, together with the Hodge filtration $F$, gives the {\it Hodge decomposition} of $L_x$ at each $x\in X$. This $F^c$ is defined on
$$\MMc:=\OO_{\XX}\otimes_{\C}L,$$
inducing a filtration of $\OO_{X_{\RRR}}\otimes L$ by scalar extension, {\it but not on} $\MM=\OO_{\XX}\otimes_{\C}\LL$ in general (unless $L$ is defined over $\RRR$). Its complex conjugate $\overline{F^c}$ is defined on
$$M^c:=\OO_X\otimes_{\C}\LL,$$
which is the complex conjugate of $\MMc$, and $\overline{F^c}$ can be described as in (A.3) below by using a {\it polarization} of variation of complex Hodge structure
$$L\otimes_{\C}\LL\to\C.
\leqno({\rm A}.2.1)$$
Note that the latter is equivalent to a $\D_X{\otimes}_{\C}\D_{\XX}$-linear pairing:
$$M\otimes_{\C}\MM\to\OO_{X_{\RRR}}\,(\subset\Db_X),
\leqno({\rm A}.2.2)$$
with $\Db_X$ is the sheaf of distributions, see also \cite[(5.1)]{ScVi}.
Here we use $\D_X{\otimes}_{\C}\D_{\XX}$-linear pairings as above. This should be distinguished from the ``sesquilinear pairing":
$$M\times M\to\OO_{X_{\RRR}}\,(\subset\Db_X),
\leqno({\rm A}.2.3)$$
which uses the {\it twist} of the action of $\OO_X$ by the complex conjugation as in (A.1.2) {\it for one factor}. We will not use sesquilinear pairings except for the case $X$ is a {\it point}, since the use of the complex conjugation in this way could be {\it very confusing} when {\it monodromies} are considered.
\msn
{\bf A.3.~Description of the opposite filtration.} We can verify that the Hodge filtration $F$ on $M=\OO_X\otimes_{\C}L$ induces the associated {\it dual} filtration $F^{\vee}$ on
$$M^c=\OO_X\otimes_{\C}\LL,$$
with respect to the polarization (A.2.1), and $F^{\vee}$ coincides, up to a shift of filtration, with the complex conjugate $\overline{F^c}$ of the opposite filtration $F^c$ on
$$\MMc=\OO_{\XX}\otimes_{\C}L.$$
\sk
Indeed, using the {\it point-wise dual filtration}, we see that the filtration $F$ on $M=\OO_X\otimes_{\C}L$ induces the filtration $F^{\vee}$ on $M^c=\OO_X\otimes_{\C}\LL$, or equivalently, $\overline{F^{\vee}}$ on $\MMc=\OO_{\XX}\otimes_{\C}L$ (by taking the complex conjugation as in (A.1.2)), and the latter coincides with the opposite filtration $F^c$ up to a shift of the filtration.
In fact, the non-degenerate pairing (A.2.1) induces an isomorphism at each $x\in X$
$$L_x\cong\LL^{\vee}_x,
\leqno({\rm A}.3.1)$$
with left-hand side having the filtration $F$, and this induces the dual filtration $F^{\vee}$ on $\LL_x$, which coincides with the opposite filtration $F^c$ up to a shift of filtration. Here it seems necessary to view the polarization as a {\it sesquilinear pairing} as in (A.2.3) after restricting it over $x\in X$, in order to get the {\it Hodge decomposition} of $L_x$ at each $x\in X$, where $L_x$ and $\LL_x$ are identified with each other up to the action of $\C$.
\msn
{\bf Remark~A.4.} It has been informed from T.~Mochizuki that a pairing corresponding to the following is used in the twistor theory (where the first factor behaves contravariantly):
$$\DD M\otimes\MMc\to\OO_{X_{\RRR}}\,(\subset\Db_X),
\leqno({\rm A}.4.1)$$
which is induced by the canonical morphism $L^{\vee}\otimes L\to\C$.
Here $\DD M$ is the dual as a filtered $\D$-module, and this is different from the dual as a filtered $\OO$-module by a shift of filtration $F$, see (A.5.2) below. Note that this shift is very important in the theory, see \cite[13.5]{Mo2}.
\msn
{\bf Remark~A.5.} In the case of variations of complex Hodge structure, it may be more natural to use the following $\D_X{\otimes}_{\C}\D_{\XX}$-linear pairing:
$$M\otimes_{\C}\overline{M'}\to\OO_{X_{\RRR}}\,(\subset\Db_X),
\leqno({\rm A}.5.1)$$
where
$$(M',F):=\DDD(M^c,\FFc[-d_X])=\DDD_{\OO}(M^c,\FFc)\,\bl(:=\HH om_{\OO_X}\bl((M^c,\FFc),(\OO_X,F)\br)\br),
\leqno({\rm A}.5.2)$$
with $\Gr_F^p\OO_X=0$ ($p\ne 0$). Indeed, the argument in (A.3) implies
$$\DDD_{\OO}(M^c,\FFc)=(M,F[w]).$$
Here $\DDD$ and $\DDD_{\OO}$ denote respectively the dual functors for filtered holonomic $\D$-modules and for filtered locally free sheaves.
\sk
Since the second factor of the pairing must behave {\it contravariantly} in the case of complex Hodge modules, the above formula (using the dual functor $\DDD$) may be useful especially in the {\it mixed} case where we have the weight filtration $W$. (Indeed, we would have to use a weight ``co-filtration" for the second factor otherwise.)
\sk
Here it seems useful to introduce the {\it Hermitian dual}\,:
$$\DDD^H(M):=\overline{\HH om_{\OO_X}(M,\Db_X)},$$
as well as the {\it Hermitian conjugate}\,:
$$M^{HC}:=\DDD(\DDD^H(M)).$$
These are isomorphic to $\DDD(M^c)$ and $M^c$ respectively in the case of variations of complex Hodge structure by the above argument. Setting $M''=\DDD(M')$, we could use the isomorphism
$$M^{HC}\cong M'',\q\h{or equivalently,}\,\,\,\DDD^H(M)\cong\DDD(M''),$$
(instead of a pairing) for the definition of complex Hodge modules, where $M''=\DDD(M')$ has the filtration $F$ (corresponding to the {\it opposite filtration} $F^c$ by the complex conjugation) so that the above isomorphism induces complex Hodge structures at any $x\in X$ in the case of polarizable variations of complex Hodge structure.
\msn
{\bf Remark~A.6.} The formulation in Remark~A.5 above seems to be useful for the proof of the $E_2$-{\it degeneration} of the weight spectral sequence, where one has to show the {\it compatibility} of the above structure with the {\it differentials} $d_r$ of the spectral sequences for any $r\ges 2$. More precisely, the {\it compatibility} of each differential $d_r$ of the two spectral sequences with the induced pairing between the $E_2$-terms (which are identified with the $E_r$-terms) must be proved {\it before showing their vanishing} for each $r\ges 2$. This problem does not seem quite easy to solve by using a {\it pairing argument} unless we adopt the above formulation together with the standard argument of spectral sequences associated with filtered complexes. Note that it is enough to prove a {\it comparison isomorphism} between two spectral sequences in the latter case. Here one would have to use the vanishing of ${\mathcal E}xt^i_{\D_X}(M,\Db_X)$ for $i\ne 0$ (see \cite{Ka}) in order to apply the derived functor argument.
\sk
As for \cite{Ka} and \cite{Sab1}, these seem to be closely related to {\it ideas} behind the theories of twistor modules and complex Hodge modules, although they do not seem to be quoted explicitly in \cite{Sab2}. (This may be rather strange, see also a sentence before \cite[(5.1)]{ScVi} or p.~19, \raise1pt\h{$\uparrow\,$}l.~6 in arXiv:1206.5547.)
\msn
{\bf Remark~A.7.} In the case of {\it real\,} Hodge modules, the above ``weight spectral sequence" is {\it defined} (or constructed) in the category $MF_{rh}(\D_X,\RRR)$ consisting of filtered regular holonomic $\D$-modules endowed with a real structure. However, the ``corresponding category" together with the ``corresponding argument" for {\it complex} Hodge modules does not seem to be very clear.
\sk
If we use the above Hermitian conjugate $M^{HC}$, then this category can be given by the category whose objects are pairs of filtered holonomic $\D$-modules $(M,F)$, $(M'',F)$ endowed with an isomorphism $\al:M^{HC}\cong M''$ as in Remark~A.5 above, and the argument is relatively easy as is explained in Remark~A.6 above.
\sk
If the Hermitian conjugate is not used, then one would have to use a {\it pairing argument} between two ``{\it spectral objects}" in the sense of Verdier \cite{Ve} in order to replace the above {\it comparison isomorphism} argument between two spectral sequences, where the argument would be much more complicated. (This pairing argument between two spectral objects seems to be useful also in the twistor case.)
\msn
{\bf Remark~A.8.} It would be highly desirable to write down {\it complete proofs} of main theorems in the theory of {\it complex} Hodge modules before giving courses about these. Otherwise, it would be rather irresponsible toward the audience.
For instance, generalizations of results of Schmid and Zucker to the {\it non-quasi-unipotent} monodromy case seem to be far from trivial.
\sk
There does not seem to be any example where rational or real Hodge modules are not enough, and complex Hodge modules are really needed in algebraic geometry, except for certain cases related to representation theory as in \cite{ScVi}.


\begin{thebibliography}{Sab4}
\bibitem[Be]{Be} Beilinson, A., How to glue perverse sheaves, L.N.M.\ 1289, Springer, Berlin (1987), 42--51.
\bibitem[BBD]{BBD} Beilinson, A., Bernstein, J.\ and Deligne, P., Faisceaux pervers, Ast\'erisque 100, S.M.F., 1982.
\bibitem[De1]{th2} Deligne, P., Th\'eorie de Hodge II, Publ.\ Math.\ IHES 40 (1971), 5--57.
\bibitem[De2]{conw} Deligne, P., La conjecture de Weil II, Publ.\ Math.\ IHES 52 (1980), 137--252.
\bibitem[Ka]{Ka} Kashiwara, M., Regular holonomic $D$-modules and distributions on complex manifolds, in Complex Analytic Singularities, Adv.\ Stud.\ in Pure Math.\ 8, North Holland, Amsterdam (1987), pp.~199--206.
\bibitem[Mo1]{Mo1} Mochizuki, T., Asymptotic behaviour of tame harmonic bundles and an application to pure twistor $D$-modules I, II. Mem.\ Am.\ Math.\ Soc.\ 185, 2007.
\bibitem[Mo2]{Mo2} Mochizuki, T., Mixed twistor $D$-Modules, L.N.M.\ 2125, Springer, Berlin, 2015.
\bibitem[Mo3]{Mo3} Mochizuki, T., Letter to the author (in Japanese), Oct.\ 19, 2016.
\bibitem[Mo4]{Mo4} Mochizuki, T., Letter to the author (in Japanese), Oct.\ 26, 2016.
\bibitem[Sab1]{Sab1} Sabbah, C., Vanishing cycles and Hermitian duality, Proc.\ Steklov Inst.\ Math.\ 238 (2002), pp.~194--214.
\bibitem[Sab2]{Sab2} Sabbah, C., Polarizable twistor $D$-modules, Ast\'erisque 300, S.M.F., 2005.
\bibitem[Sab3]{Sab3} Sabbah, C., Isomonodromic deformations and Frobenius manifolds - An introduction, Universitext, Springer, London, 2007.
\bibitem[Sab4]{Sab4} Sabbah, C., Exponential Hodge theory, arXiv:1511.00176.
\bibitem[Sai1]{bl} Saito, M., On the structure of Brieskorn lattice, Ann.\ Inst.\ Fourier 39 (1989), 27--72.
\bibitem[Sai2]{mhm} Saito, M., Mixed Hodge modules, Publ. RIMS, Kyoto Univ.\ 26 (1990), 221--333.
\bibitem[Sai3]{ams} Saito, M., Arithmetic mixed sheaves, Inv.\ Math.\ 144 (2001), 533--569.
\bibitem[Sai4]{tdef} Saito, M., Twistor deformation of rank one local systems, arXiv:1307.4907.
\bibitem[Sai5]{ypg} Saito, M., A young person's guide to mixed Hodge modules, arXiv:1605.00435.
\bibitem[ScSt]{ScSt} Scherk, J.\ and Steenbrink, J.H.M., On the mixed Hodge structure on the cohomology of the Milnor fibre, Math.\ Ann.\ 271 (1985), 641--665.
\bibitem[ScVi]{ScVi} Schmid, W.\ and Vilonen, K., Hodge theory and unitary representations of reductive Lie groups, in Frontiers of mathematical sciences, Int.\ Press, Somerville, MA, 2011, pp.~397--420 (or arXiv:1206.5547).
\bibitem[Si]{Si} Simpson, C., Mixed twistor structures, arXiv:math.AG/9705006.
\bibitem[St]{St} Steenbrink, J.H.M., Limits of Hodge structures, Inv.\ Math.\ 31 (1976), 229--257.
\bibitem[Va]{Va} Varchenko, A.N, Asymptotic Hodge structure in the vanishing cohomology, Math.\ USSR-Izv.\ 18 (1982), 469--512.
\bibitem[Ve]{Ve} Verdier, J.-L., Des cat\'egories d\'eriv\'ees des cat\'egories ab\'eliennes, Ast\'erisque 239 (1996).
\end{thebibliography}
\end{document}